\documentclass[a4paper, 11pt, abstract, DIV12]{scrartcl}

\makeatletter
\DeclareOldFontCommand{\rm}{\normalfont\rmfamily}{\mathrm}
\DeclareOldFontCommand{\sf}{\normalfont\sffamily}{\mathsf}
\DeclareOldFontCommand{\tt}{\normalfont\ttfamily}{\mathtt}
\DeclareOldFontCommand{\bf}{\normalfont\bfseries}{\mathbf}
\DeclareOldFontCommand{\it}{\normalfont\itshape}{\mathit}
\DeclareOldFontCommand{\sl}{\normalfont\slshape}{\@nomath\sl}
\DeclareOldFontCommand{\sc}{\normalfont\scshape}{\@nomath\sc}
\makeatother

\usepackage[utf8]{inputenc}
\usepackage[T1]{fontenc}
\usepackage{lmodern}
\usepackage[draft=false, kerning=true]{microtype}

\usepackage[dvipsnames]{xcolor}
\usepackage[english]{babel}
\usepackage{amsmath, amsfonts, amsthm}
\usepackage{bm}
\usepackage{mathtools}
\usepackage[shortlabels]{enumitem}
\usepackage{physics}
\usepackage{enumitem}
\usepackage{array, makecell}
\usepackage{multicol}
\usepackage{tablefootnote}
\usepackage{multirow, makecell}

\usepackage{tikz}
\usetikzlibrary{calc}
\newcommand{\tikzmark}[1]{\tikz[overlay,remember picture] \node (#1) {};}
\newcommand{\DrawBox}[3][]{%
  \tikz[overlay,remember picture]{
    \draw[#1]
    ($(#2)+(-0.20em,0.55em)$) rectangle
    ($(#3)+(0.25em,-0.1em)$);}
}
\newcommand{\DrawLine}[3][]{%
  \tikz[overlay,remember picture]{
    \draw[#1]
    ($(#2)+(0.185em,0.55em)$) --
    ($(#3)+(0.185em,-0.1em)$);}
}
\newcommand{\DrawLineCorr}[3][]{%
  \tikz[overlay,remember picture]{
    \draw[#1]
    ($(#2)+(0.40em,0.55em)$) --
    ($(#3)+(0.185em,-0.1em)$);}
}
\newcommand{\DrawLineCorrTwo}[3][]{%
  \tikz[overlay,remember picture]{
    \draw[#1]
    ($(#2)+(0.50em,0.55em)$) --
    ($(#3)+(0.185em,-0.1em)$);}
}

\newcommand{\N}{\mathbb{N}}
\newcommand{\Z}{\mathbb{Z}}
\newcommand{\F}{\mathbb{F}}

\DeclareMathOperator{\ag}{AG}
\DeclareMathOperator{\pg}{PG}
\DeclarePairedDelimiter\set\{\}
\DeclarePairedDelimiterX{\setm}[2]{\{}{\}}{#1\,\delimsize\vert\,\mathopen{}#2}
\let\abs\relax
\DeclarePairedDelimiter{\abs}{\lvert}{\rvert}

\newtheorem{theorem}{Theorem}
\newtheorem{definition}{Definition}[section]

\newtheorem{corollary}[theorem]{Corollary}

\theoremstyle{remark}

\newtheorem{remark}[definition]{Remark}

\numberwithin{equation}{section}

\makeatletter
\def\blfootnote{\xdef\@thefnmark{}\@footnotetext}
\makeatother

\usepackage{etoolbox}
\newcounter{savefootnote}\newcounter{savempfootnote}
\newcounter{symfootnote}\newcounter{symmpfootnote}
\AtBeginEnvironment{minipage}{\setcounter{symmpfootnote}{0}}
\newcommand{\symfootnote}[1]{%
   \setcounter{savefootnote}{\value{footnote}}\setcounter{savempfootnote}{\value{mpfootnote}}%
   \setcounter{footnote}{\value{symfootnote}}\setcounter{mpfootnote}{\value{symmpfootnote}}%
   \ifnum\value{footnote}>8\setcounter{footnote}{0}\fi\ifnum\value{mpfootnote}>8\setcounter{mpfootnote}{0}\fi%
   \let\oldthefootnote=\thefootnote\let\oldthempfootnote=\thempfootnote%
   \renewcommand{\thefootnote}{\fnsymbol{footnote}}\renewcommand{\thempfootnote}{\fnsymbol{mpfootnote}}%
   \footnote{#1}%
   \let\thefootnote=\oldthefootnote\let\thempfootnote=\oldthempfootnote%
   \setcounter{symfootnote}{\value{footnote}}\setcounter{symmpfootnote}{\value{mpfootnote}}%
   \setcounter{footnote}{\value{savefootnote}}\setcounter{mpfootnote}{\value{savempfootnote}}%
}

\usepackage{hyperref}
\hypersetup{
  bookmarksopen=true,
  bookmarksnumbered=true,
  pdftitle={Exponentially Larger Affine and Projective Caps},
  pdfauthor={Christian Elsholtz, Gabriel F. Lipnik},
  pdflang=english,
  colorlinks = false,
  linkbordercolor = {MidnightBlue},
  urlbordercolor = {MidnightBlue},
  citecolor = {ForestGreen},
  citebordercolor = {ForestGreen!80!black},
  plainpages=false
}
\usepackage{xurl}
\hypersetup{breaklinks=true}

\begin{document}

\title{\vspace*{-2ex}
  Exponentially Larger\\
  Affine and Projective Caps%
  \footnotetext{%
  \begin{description}
  \item [Christian Elsholtz] \texttt{elsholtz@math.tugraz.at},
    Institute of Analysis and Number Theory, Kopernikusgasse~24/II,
    Graz University of Technology, 8010 Graz, Austria
  \item [Gabriel F.\ Lipnik]
    \href{mailto:math@gabriellipnik.at}{\texttt{math@gabriellipnik.at}},
    Institute of Analysis and Number Theory, Kopernikusgasse~24/II,
    Graz University of Technology, 8010 Graz, Austria
  \item [Support] The authors acknowledge the support of the Austrian
    Science Fund (FWF): W\,1230 and I\,4945-N. The first author is
    also supported by grant HR 12/2020, WTZ Kroatien S\&T Croatia
    2020\_21, awarded by the OeAD.
  \item [Acknowledgement] The authors would like to thank Seva Lev for
    sending his unpublished manuscript~\cite{lev}.
  \item [2020 Mathematics Subject Classification]  51E20, 51E22, 05B25
  \item[Key words] caps, affine space, projective space,
    progression-free sets
  \end{description}}}

\author{Christian Elsholtz and Gabriel F.\ Lipnik}
\date{\vspace{-3ex}}
\maketitle

\begin{abstract}
In spite of a recent breakthrough on upper bounds of the size of cap sets
(by Croot, Lev and Pach (2017) and Ellenberg and Gijswijt (2017)),
the classical cap set constructions had not been affected. In this work, we introduce a very different method of construction for caps in all affine spaces with odd prime modulus~$p$.  Moreover, we show that for all 
primes $p \equiv 5 \bmod 6$ with $p \leq 41$, the new construction leads to an exponentially larger growth of the affine and projective caps in $\ag(n,p)$ and $\pg(n,p)$. For example, when $p=23$, the existence of caps with growth $(8.0875\ldots)^n$ follows from a three-dimensional example of Bose (1947), and the only improvement had been to $(8.0901\ldots)^n$ by Edel (2004), based on a six-dimensional example. We improve this lower bound to $(9-o(1))^n$.
\end{abstract}

\section{Introduction and Overview}
The study of large point sets without three points on any line, in affine or projective spaces, is a classical topic in geometry, and more recently also in additive combinatorics. 
An introduction and some general information on these sets called \emph{caps}, in particular from a geometric point of view, can be found in several chapters of Hirschfeld's three volumes on projective geometries over finite fields~\cite{hirschfeld-1985-finite-projective-spaces, hirschfeld-1998-finite-geometries, hirschfeld-thas-2016-general-galois-geometries}, in a survey by Hirschfeld and Storme~\cite{Hirschfeld-Storme}, in relevant papers by Bierbrauer and Edel,
e.g.~\cite{bierbrauer-2003-large-caps, bierbrauer-edel-2014-large-caps-in-projective-galois-spaces, edel:2004:product-caps, edel-bierbrauer-1999-recursive-constructions-for-large-caps},
and on Edel's website~\cite{edel-website}.

A lot of results study the size of complete caps (i.e., caps which cannot be extended) in a fixed dimension over a fixed finite field; see e.g.~\cite{MR3093866,MR3684305, MR3800841}.
It is even an open problem to characterize complete caps in dimension~3 over~$\F_q$; see for example Hirschfeld and Thas~\cite{Hirschfeld-Thas}. Numerous papers give alternative constructions for non-equivalent caps; see e.g.\ Kroll and Vincenti~\cite{kroll-vincenti}.

An important breakthrough 
\cite{bennett, croot-lev-pach:2017:progression-free-sets-in-Z4, ellenberg-gijswijt:2017:subsets-without-3-term-arithmetic-progression, grochow}
has recently lead
to greatly improved upper bounds for the largest possible size of these sets in the affine geometry $\ag(n,p)$. 

In this paper, we improve longstanding \emph{lower bounds} for caps when $p\in \{11,17,23,29,41\}$.
In fact, the improvement is actually an exponential improvement (in the standard terminology, see~\cite{croot-lev-pach:2017:progression-free-sets-in-Z4}). It might be clear that this does not come from a refinement of previous methods but from an entirely different approach.

Previous cap set constructions are (recursively) based on a product construction from good examples in low dimensions, which we think of as a ``local'' approach; see for example~\cite{bierbrauer-2003-large-caps, bierbrauer-edel-2014-large-caps-in-projective-galois-spaces,
Davydov-Ostergard,
edel:2004:product-caps,
edel-bierbrauer-1999-recursive-constructions-for-large-caps, 
Mukhopadhyay}.
In contrast, we construct
a set of vectors with certain constraints with regards to the occurring digits, similar to a construction by Salem and Spencer~\cite{salem-spencer:1942:sets-without-3-term-arithmetic-progression} in the integer case, and we think of this construction as a ``global'' approach. 


In this paper, we describe a new type of cap construction
in the affine space $\ag(n,p)$ over the field $\Z_p$ with $p \geq 5$ prime (and therefore also in the corresponding projective space $\pg(n,p)$)
that actually works for all dimensions
over $\Z_p$. In its most basic case
this includes the simple cap construction $\{0,1\}^n \subset \F_3^n$. This has been generalized previously to certain product constructions. In this paper, we generalize this in a novel way to combine well-chosen digit sets with certain conditions.
It will be apparent from the construction below that for given $n$ and $p$, there are usually many non-equivalent caps; see Section~\ref{sec:non-equivalent}. For some primes, we can even achieve new records of the largest known caps and we will concentrate on this aspect. 
These appear to be the first improvements over the results of Bose, Bierbrauer and Edel; for details see below.

In the following, we consider the affine space $\ag(n,p)$, where~$n\in\N$ is the dimension and~$p$ is a prime (and thus, $\ag(n,p) = ((\Z/p\Z)^n,+)$, or $\Z_p^n$ for brevity).
An affine cap~$S$ is a subset of~$\Z_p^n$ such that no three points in~$S$ are collinear, i.e., for any three pairwise distinct points  $x$, $y$, $z\in S$, the vectors $y - x$ and $z - x$ are linearly independent over~$\Z_p$.
This condition is equivalent to the fact that for any $(a,b,c)\in \Z_p^3\setminus\{(0,0,0)\}$ with 
$a+b+c=0$, one also has $ax + by + cz \neq 0$.

Projective caps are analogue sets in the projective space~$\pg(n,p)$ instead of~$\ag(n,p)$. Since affine spaces can be embedded into projective spaces, our improved caps in~$\ag(n,p)$ also represent caps in the projective space.

\paragraph{Related Work.} It is known that for $m\in \{3,4,5\}$, a cap in $\Z_m^n$ is
equivalent to a set in which no three distinct points are in an
arithmetic progression.  (Note that $\Z_4^n \neq \ag(n,p)$.) There were
important contributions by Brown and
Buhler~\cite{brown-buhler:1982:density-geometric-ramsey}, Frankl,
Graham and
Rödl~\cite{frankl-graham-roedl:1987:no-3-term-arithmetic-progression},
Meshulam~\cite{meshulam:1995:no-3-term-arithmetic-progression},
Lev~\cite{lev:2004:progression-free-sets}, Bateman and
Katz~\cite{bateman-katz:2012:bounds-on-cap-sets}, Croot, Lev and
Pach~\cite{croot-lev-pach:2017:progression-free-sets-in-Z4},
Ellenberg and
Gijswijt~\cite{ellenberg-gijswijt:2017:subsets-without-3-term-arithmetic-progression} as well as Petrov and Pohoata~\cite{petrov-pohoata-2020-improved-bounds-for-progression-free-sets}. Moreover,
some readers may recall the case of caps in~$\Z_{3}^{4}$ from the
popular card game \emph{SET}~\cite{davis-maclagan:2003:set}.

So far, the best-known approach to construct caps \emph{for general prime modulus} is
to take a simple product construction of large caps in low dimension. 
Let~$C_{n,p}$ and $C_{n,p}^{\rm{pr}}$ denote\footnote{Note that for projective caps the size of the largest cap
is often denoted by $m_2(n,p)$, and sometimes $m_2^{\rm{aff}}(n,p)$ is used in the affine case.} the sizes of the largest affine and projective cap in dimension~$n$, respectively.
It is known that the largest affine cap in dimension~$3$ has size~$p^2$, i.e., $C_{3, p} = p^2$; see for example in Bierbrauer~\cite{bierbrauer-2003-large-caps}.
These maximal caps are also called ovaloids.
In PG$(3,q)$ with odd $q$, these maximal caps come from elliptic quadrics, see \cite{O'Keefe}.
A representative for such a cap is the set
\[ \setm{(t^2+st+as^2,s,t,1)}{s,t \in \F_q} \cup \{(1,0,0,0) \},\]
where $x^2+x+a$ is irreducible over $\F_q$.
In the corresponding affine space, the point $(1,0,0,0)$ is removed.
 As a consequence, we obtain the bound $C_{n,p} \gg p^{2n/3}$ by simply taking products of this cap.
This result can be considered classical,
as the determination of the size of caps in $\pg(3,q)$ for odd prime powers $q$ goes back to Bose~\cite{bose} in 1947. 

The refinement by Edel and Bierbrauer is based on the fact
that one can form an almost-product of special projective caps, namely if they possess a tangent hyperplane (see~\cite[Theorem~10]{edel-bierbrauer-1999-recursive-constructions-for-large-caps}). In particular, this gives
$(q^2+1)^2-1=q^4+2q^2$ points in PG$(6,q)$, and the reduction to the affine space gives $q^4+q^2-1$ points in AG$(6,q)$; see~\cite[Section 1]{edel:2004:product-caps}.

There are several computational results on caps in small dimension; see \cite{edel-website,Hirschfeld-Storme}. However, the only 
known asymptotic improvement over Bose's result on the lower bound when $p \geq 5$ is due to Bierbrauer and 
Edel~\cite[Theorem~11]{edel-bierbrauer-1999-recursive-constructions-for-large-caps}
for projective caps 
and Edel~\cite{edel:2004:product-caps} for affine caps, and is based on a product construction of a large cap in dimension~$6$. If~$n$ is a multiple of~$6$, then
 Edel's construction yields $C_{n,p}^{\rm{pr}}\geq C_{n,p}\geq (p^4+p^2-1)^{n/6}$. If~$n$ is not a multiple of~$6$,
 one can modify the construction slightly, but in any case this only
 influences a constant $C_{p}$ in $C_{n,p}\geq C_p (p^4+p^2-1)^{\lfloor n/6\rfloor}$.
 
It is known that the limit $c_{p}=\lim_{n \rightarrow \infty} (C_{n,p})^{1/n}$ exists and is in the interval $[2,p)$; see for example~\cite[Proposition~3.8]{elsholtz-pach:2019:progression-free-sets-caps}.
Numerically, Edel's construction gives 
only a small improvement of the earlier bound $p^{2n/3}$. For example, when
$p=17$, then the bound
$6.611\ldots$ is improved to $6.615\ldots\,$. In this paper, we will improve this to~$7$.
When $p=23$, a lift of Bose's result gives
a constant $8.087\ldots\,$, which Edel improved to $8.090\ldots\,$.
We improve this to~$9$. However, while Edel's construction works for all primes, our construction has to be optimized for individual primes.

For $p=5$, Edel's construction gives $c_5\geq 649^{1/6}=2.942\ldots$\,. 
Recently, Elsholtz and Pach~\cite{elsholtz-pach:2019:progression-free-sets-caps} have constructed large progression-free sets and
it emerged that in~$\Z_5$, their construction is asymptotically better than Edel's bound; Edel's lower bound was improved to
$c_{5}\geq 3$. In the case modulo~$4$ (i.e., working in $\Z_4^n$ rather than
$\F_4^n$), Elsholtz and Pach~\cite{elsholtz-pach:2019:progression-free-sets-caps} gave a much more substantial improvement from 
$c_{4}\geq 2.519\ldots$ to $c_{4}\geq 3$. Improvements in the case of a prime base
seem to be much more difficult, since the existing construction of Edel seems to be good.

Another important measure for the size of caps is the 
exponent~$\mu(p)=\lim_{n \rightarrow \infty} 
(\log_p C_{n,p})/n$  in the representation of the size 
as~$p^{\mu(p) n}$. The mentioned result $C_{n,p} \gg p^{2n/3}$ clearly implies $\mu(p) \geq 2/3$. 
The recent breakthrough of Ellenberg and Gijswijt~\cite{ellenberg-gijswijt:2017:subsets-without-3-term-arithmetic-progression}
shows that $\mu(p)<1$.
Indeed, their method yields the bound
\begin{equation*}
    C_{n,p}\leq (J(p)p)^n,
\end{equation*}
where 
\begin{equation*}
J(p)=\frac{1}{p}\min\limits_{0<t<1} \frac{1-t^p}{(1-t)\, t^{(p-1)/3}};
\end{equation*}
see~\cite{Blasiak-Church-etal}. It is known that $J(p)$ is decreasing and $\lim_{p \rightarrow \infty} J(p) =0.8414\ldots\,$; see~\cite[Equation~(4.11)]{Blasiak-Church-etal}.

Besides the mentioned product constructions, also another approach is known: In an unpublished work of Lev~\cite{lev}, he describes an elegant method to ``globally'' construct large caps in $\F_3^n$. These caps have basically the form
\begin{equation*}
    S = \setm{(x,y,x^2 - \lambda y^2)}{x, y\in \F_q} \subseteq \F_q^3 \cong \F_3^n,
\end{equation*}
where $\lambda\in\F_q$ is a fixed non-square, $n$ is a multiple of~$3$ and $q = 3^{n/3}$. However, these sets have size $3^{2n/3} = (2.08008\ldots)^n$, which is of the same quality as Bose's construction lifted to higher dimension.

\paragraph{Overview of our Work.} In this paper, we extend the combinatorial method 
of Elsholtz and Pach~\cite{elsholtz-pach:2019:progression-free-sets-caps} from the case of sets avoiding arithmetic progressions to affine caps (with prime modulus larger than~$4$). In particular, we introduce some new directions for finding good digit sets, which are crucial for our constructions of large caps; see Section~\ref{sec:admissible-sets}.

Our results improve the lower bounds of $c_{p}$ for
$p\in\set{11,17,23,29,41}$, and the improvements in these cases are indeed
substantial. Especially the case $p=23$ with an exponent of $\mu(23) \geq 0.70075\ldots$ comes quite close to the case of $p=3$, where a construction is known based on a large cap in dimension~480, giving $\mu(3)\geq 0.72485\ldots$\,; see~\cite[Section~5]{edel:2004:product-caps}.

Table~\ref{tab:comparison} compares our new lower bounds
to those by Edel \cite{edel:2004:product-caps}.

\begin{table}[htbp]
\begin{minipage}{\textwidth}
  \centering
    \begin{tabular}{r|r|r|r|r|r}
    \makecell{\multirow{2}{*}{\thead{$\bm{p}$}}}& \multicolumn{4}{c|}{\thead{\textbf{lower bounds for~$\bm{c_p}$}}} & \makecell{\multirow{2}{*}{\thead{\textbf{exponent} $\bm{\mu(p)}$}}}\\
      &\thead{$p^{2/3}$} & \thead{$(p^4 + p^2 - 1)^{1/6}$} & \makecell{\thead{new}} & \makecell{\thead{improvement}\symfootnote{Compared to the best previously known bound $(p^4 + p^2 - 1)^{1/6}$.}}\\\hline
      $5$  & $2.92401\ldots$ & $2.94243\ldots$ & $\bm{3}$
      & $1.9562\%$ & $\bm{0.68260\ldots}$\\
      $7$ & $3.65930\ldots$ & $3.67139\ldots$ & $3$ && $0.56457\ldots$\\
      $11$ & $4.94608\ldots$ & $4.95282\ldots$ & $\bm{5}$ & $0.9526\%$ & $\bm{0.67118\ldots}$\\
      $13$ & $5.52877\ldots$ & $5.53418\ldots$ & $4$ & & $0.54047\ldots$\\
      $17$ & $6.61148\ldots$ & $6.61528\ldots$ & $\bm{7}$ & $5.8156\%$ & $\bm{0.68682\ldots}$\\
      $19$ & $7.12036\ldots$ & $7.12364\ldots$ & $6$ & & $0.60852\ldots$\\
      $23$ & $8.08757\ldots$ & $8.09012\ldots$ & $\bm{9}$ & $11.2468\%$ & $\bm{0.70075\ldots}$\\
      $29$ & $9.43913\ldots$ & $9.44099\ldots$ & $\geq\bm{10}$ & $\geq 5.9210\%$ & $\geq \bm{0.68380\ldots}$\\
      $31$ & $9.86827\ldots$ & $9.86998\ldots$ & $\geq 8$ & & $\geq 0.60554\ldots$\\
      $37$ & $11.10370\ldots$ & $11.10505\ldots$ & $\geq 10$ & & $\geq 0.63767\ldots$\\
      $41$ & $11.89020\ldots$ & $11.89138\ldots$ & $\geq\bm{12}$ & $\geq 0.9134\%$ & $\geq \bm{0.66914\ldots}$\\
    \end{tabular}
    \end{minipage}
  \caption{Comparison of previously known best lower bounds for $c_p$ to our new ones, and new lower bounds for the exponent~$\mu(p)$. 
  Figures in bold constitute new records.
  The $\geq$-sign is meant to indicate cases in which we cannot ensure that our method is not able to produce better results than the stated ones.}
  \label{tab:comparison}
\end{table}

\section{Results and Construction}
\label{sec:result}

In the following, we use Vinogradov's notation, 
where  $f(n)\gg g(n)$ means that there exists some $C>0$ such that
$f(n) \geq C g(n)$ holds for all $n > n_0$.

We directly start by stating our main result.

\begin{theorem}
  \label{thm:main}
  If $C_{n,p}$ denotes the size of the largest affine cap in~$\Z_p^n$
  and $c_{p} = \lim_{n\to\infty}(C_{n,p})^{1/n}$, then the following holds:
  \begin{equation*}
    C_{n,11}\gg \frac{5^{n}}{n^{1.5}},\quad C_{n,17} \gg \frac{7^{n}}{n^{2.5}},\quad C_{n,23} \gg
    \frac{9^{n}}{n^{3.5}},\quad C_{n,29}\gg \frac{10^{n}}{n^{4}}\qq{and}
    C_{n,41}\gg \frac{12^{n}}{n^{5}}.
  \end{equation*}
  As a consequence, we have that
  \begin{equation*}
    c_{11} \geq 5,\quad c_{17} \geq 7,\quad c_{23}\geq 9,\quad c_{29}\geq 10
    \qq{and} c_{41}\geq 12.
  \end{equation*}
  \end{theorem}

Since every subset of~$\ag(n,p)$ can be embedded into~$\pg(n,p)$, this directly implies the following corollary.
\begin{corollary}
 The lower bounds from Theorem~\ref{thm:main} also hold for the largest caps in~$\pg(n,p)$.
\end{corollary}
Moreover, our improved bounds on caps can also be transformed into improved bounds on linear codes. For details we refer to~\cite[Theorem~1]{edel-bierbrauer-1999-recursive-constructions-for-large-caps}.

These new bounds are based on a ``global'' construction of affine caps: We take a set of
$n$-dimensional points, where the set depends on $n$ in a much stronger way
than taking a tensor product construction of a small (local) cap. The idea
is, for a fixed prime~$p$, to find a large set of digits $D\subseteq \Z_{p}$ and a subset $D'\subseteq D$ such that the set 
\begin{equation}
\label{eq:def-of-S}
  S(D, D', n) \coloneqq \setm[\Bigg]{(a_{1},\ldots,a_{n})\in D^{n}}{\forall d\in D'\colon a_{i}
    = d \text{ for } \frac{n}{\abs{D}} \text{ values of $i$}}
\end{equation}
is a cap in $\ag(n,p)$ for all $n\in\N$ with $\abs{D}\mid n$. If this is the case, then we say that $(D, D')$
is \emph{admissible}. Moreover, we say that~$D$ is admissible if there is some $D'\subseteq D$ such that $(D, D')$ is admissible. Note that if $D_1' \subseteq D_2' \subseteq D$, then the admissibility of $(D, D_1')$ implies the admissibility of~$(D, D_2')$. 

Next, we combinatorially determine the cardinality of the set $S(D, D', n)$ and then asymptotically estimate it by applying Stirling's formula, which leads to
\begin{equation}
\label{eq:cap-size}
    \abs{S(D, D', n)} = \Biggl(\prod_{\ell=0}^{\abs{D'}-1}\binom{n - \frac{\ell n}{\abs{D}}}{\frac{n}{\abs{D}}}\Biggr)(\abs{D} - \abs{D'})^{n-\frac{\abs{D'}n}{\abs{D}}}  \sim \frac{c\abs{D}^{n}}{n^{\delta/2}}
\end{equation}
with 
\begin{equation*}
     \delta = \min\set{\abs{D'}, \abs{D} - 1} \qq{and} c = \frac{1}{\sqrt{1 - \delta/\abs{D}}} \biggl(\frac{\abs{D}}{2\pi}\biggr)^{\delta/2}.
\end{equation*}
The form of the parameter~$\delta$ comes from the fact that fixing the frequencies of~$\abs{D}$ digits leads to the same result as fixing the frequencies of~$\abs{D} - 1$ digits, because then the frequency of the last digit is fixed automatically.  With the usual interpretation $0^0 = 1$, \eqref{eq:cap-size} also holds true for $D' = D$. The given cardinality is of order $(\abs{D}-o(1))^n$ as~$n$ increases.

In order to obtain a large cap, \eqref{eq:cap-size} implies that, first of all, we need to
\begin{itemize}
    \item choose the digit set~$D$ as \emph{large} as possible, and then
    \item find a corresponding set~$D'\subseteq D$ of digits with fixed frequencies which is as \emph{small} as possible.
\end{itemize}
However, the minimization of the set~$D'$ is restricted by the fact that the frequency conditions are crucial to ensure that the resulting set is indeed a cap. More details can be found in Section~\ref{sec:admissible-sets}.

Finally, we give some additional comments on the construction.

\begin{remark}
  \begin{enumerate}[(a)]
    \item For simplicity, the reader can assume in the first reading that $D' = D$. This still covers all the main improvements, and only slightly weakens the exponent of~$n$ in the denominators of our results.
    \item It is not crucial for our method that the frequencies of the digits in~\eqref{eq:def-of-S} are exactly $n/\abs{D}$. Other constants which can also vary depending on the digit and add up to $n$ also work. However, if we want to maximize the size of the cap $S(D, D', n)$, then $n/\abs{D}$ is the best choice, in view of the multinomial distribution.
      \item If the dimension~$n$ is not a multiple of~$\abs{D}$, then we can trivially extend the set~$S(D,D', n - (n\bmod \abs{D}))$ to a subset of~$\Z_p^n$ by filling the remaining coordinates with a good cap in dimension~$n\bmod \abs{D}$. As a consequence, \eqref{eq:cap-size} holds for all~$n\in\N$, understood as an asymptotic lower bound with a slightly weaker constant~$c$.
      \item One could also think of restrictions other than fixing the frequency of some digits, e.g., fixing the ``radius'' of the points (compare Behrend's construction for progression-free sets, the application to the multidimensional setting as explained by Petrov and Poahata~\cite{petrov-pohoata-2020-improved-bounds-for-progression-free-sets} in the case modulo~$8$
    and by Elsholtz and Pach
      \cite{elsholtz-pach:2019:progression-free-sets-caps} more generally). Or one could think of fixing the frequency of multiple digits \emph{together} (as mentioned in~\cite[Proof of Theorem~3.11]{elsholtz-pach:2019:progression-free-sets-caps}). Both approaches do not seem to work for caps \emph{in general}. However, we have refrained from further optimizing the denominators in Theorem~\ref{thm:main}.
      \item It turned out that if $D$ is admissible, then the corresponding set $D'$ can be chosen in such a way that $\abs{D'} \leq \abs{D} - 2$ holds. We believe that this is always possible.
      \item So far, our method only leads to an improvement for small primes~$p$ with
  $p\equiv 5 \bmod 6$. It would be nice to have improved constructions for many primes.
  
  \item It seems to be possible to add some smaller caps to a large cap constructed in this way so that the union of all points is still a cap. This would improve the constant~$c$ by a small factor (probably less than $2$). For some details see~\cite[Theorem~3.2 and Corollary~3.4]{elsholtz-pach:2019:progression-free-sets-caps}.
  \end{enumerate}
\end{remark}

\section{Approaches for Finding Admissible Sets}
\label{sec:admissible-sets}

As already mentioned in the introduction, for $p = m\in\set{3,4,5}$
the cap set condition can be verified by only ensuring that no three points $x$,
$y$ and $z$ from the set satisfy $x + z = 2y$ (which describes
arithmetic progressions).  For $p > 5$, the cap
set condition is not only based on this equation, but also on the other equations
$ax+by+cz=0$, where $a$,
$b$, $c\in\Z_{p}$
with $a+b+c=0$.  If $m=p$ is a prime, then without loss of generality, it
is enough to assume that $a=1$. With
$c=-(b+1)$ we can assume that
$b\in \{1, \ldots, p-2\}$. (If $b=0$, then
$c=-1$ simply means that~$x$, $y$ and~$z$ are distinct. If
$b=p-1$, then we have $c=0$ with the same
consequence.)

\subsection{Modelling the Problem}

For the moment, let $b\in\set{1,\ldots,p-2}$ be fixed and
$c = - (b+1)$.  Moreover, let
\begin{equation*}
  P_{b}(D) = \setm{(x,y,z) \in D^{3}}{x + by + cz = 0 \text{ and not } x=y=z}
\end{equation*}
be the set of non-trivial ``weighted progressions'' corresponding
to~$b$. Assume that there is some $n\in\N$ with $\abs{D}\mid n$ such that there are three points $x= (x_1,\ldots,x_n)^\top$, $y=(y_1,\ldots,y_n)^\top$, $z=(z_1,\ldots,z_n)^\top\in S(D, D', n)$ which lie on a
line. For each weighted progression~$v = (v_{1},v_{2},v_{3})\in P_{b}(D)$, we
introduce a variable~$\chi_{v}$ which describes the number of occurrences of~$v$ in
the components of these three points, i.e.,
\begin{equation*}
    \chi_v = \abs[big]{\setm[\big]{i\in \set{1,\ldots,n}}{(x_i,y_i,z_i) = v}}.
\end{equation*}
Because every
digit~$d$ in~$D'$ has to occur the same number of times, we find the
equations
\begin{equation}
  \label{eq:equations-for-progressions}
  \sum_{\substack{v\in P_{b}(D)\\ v_{1} = d}}\chi_{v} = \sum_{\substack{v\in P_{b}(D)\\ v_{2} = d}}\chi_{v} \qq{and} \sum_{\substack{v\in P_{b}(D)\\ v_{1} = d}}\chi_{v} = \sum_{\substack{v\in P_{b}(D)\\ v_{3} = d}}\chi_{v} 
\end{equation}
for each~$d\in D'$.

Now it is easy to see that the non-existence
of a non-negative non-trivial integral solution
$\chi = (\chi_{v}\mid v\in P_{b}(D))$ for the equations above for all $b\in\set{1,\ldots,p-2}$ is equivalent to
the non-existence of three points on a line, i.e., the fact that
$S(D, D', n)$ is a cap.
So in order to prove the admissibility of some~$(D, D')$, we have to
ensure that the polyhedron
\begin{equation*}
  \mathcal{P} = \setm[\big]{\chi\in\Z_{\geq 0}^{\ell}}{A\cdot\chi=0}
\end{equation*}
only contains the zero vector for all $b\in\set{1,\ldots,p-2}$, where the system of linear
equations~$Ax = 0$ describes the equations given
in~\eqref{eq:equations-for-progressions} and clearly depends on~$b$ and $\ell = \abs{P_b(D)}$. This can be done by methods
of linear integer programming, e.g. with a standard IP solver. For this article, we have used the MILP packages of SageMath~\cite{SageMath:2020:9.0} as well as JuMP, an optimization package of Julia~\cite{bezanson-2017-julia}. A complete list of all admissible digit sets of maximal size for small~$p$ can be found at \url{https://gitlab.com/galipnik/large-caps}.

One way of ensuring that an admissible digit set has largest size, say size~$\ell$, among all admissible digit sets for fixed~$p$ is to find a feasible solution of the IP for \emph{all}\footnote{However, some equivalent digit sets can be neglected. For example, we can always assume without loss of generality that an admissible digit set contains the digits~0 and~1.} possible digit sets of size~$\ell + 1$ for at least one~$b$ (which implies that all these sets cannot be admissible). We have done this for $p\leq 23$; see also Table~\ref{tab:comparison}. In order to give an idea of the computation times, our implementation took about 95~minutes for the case $p = 23$ (and $\ell + 1 = 10$), while it was executed on an Intel(R) Core(TM) i7-7500U CPU at 2.70GHz. In other words, showing the non-admissibility of thousands of individual digit patterns each took only a fraction of a second.

Unfortunately, deciding if a
polyhedron contains an integer point is NP-complete~\cite{garey-johnson-computers-and-intractability}, which implies, together with the fact that the number of possible digit sets also grows exponentially for increasing~$p$,
that checking admissibility for all possible digit sets modulo~$p$
can only be done for small~$p$. Hence, it is very natural to look for simpler ways of checking whether digit sets are admissible. Two such approaches are described in the following section.

For an illustration of setting up the equations given in~\eqref{eq:equations-for-progressions} as well as the corresponding constraint matrix~$A$, we refer to the case $p = 23$ in Section~\ref{sec:proof}.
  
\subsection{Digit-Reducibility as a Sufficient Condition}
\label{sec:digit-red}

Besides the computational method presented in the previous section, we
next give a sufficient condition for the admissibility of a digit set,
which allows us to verify very easily that a set is admissible.

A pair~$(D, D')$ with $D'\subseteq D\subseteq\Z_p$ is said to be \emph{digit-reducible} if for every
$b\in\set{1,\ldots,p-2}$ and $c = -(b + 1)$ the following recursively defined algorithm
results in the empty set: If there
exists a position $r \in \{1,2,3\}$ and a digit~$d\in D'$ such that~$d$ does not
occur at position~$r$ in any of the triples in~$P_{b}(D)$ but it occurs at one of the other positions in at least one of the triples in~$P_{b}(D)$, then remove all
weighted progressions from~$P_{b}(D)$ which contain~$d$ at \emph{any position}.
Recursively apply this
rule to the remaining set~$P_{b}(D)$ again. If there do not exist an
$r \in \{1,2,3\}$ and a digit~$d\in D'$ such that~$d$ does not occur in any of
the triples in~$P_{b}(D)$ at position~$r$ but it occurs in at least one triple at any position, then stop the process.

We now explain why the reducibility of~$(D, D')$ implies that $S(D, D', n)$ is a cap for all $n\in\N$ with $\abs{D}\mid n$. Assume that~$(D, D')$ is reducible and there are three pairwise different
vectors $x = (x_1,\ldots,x_n)^\top$, $y=(y_1, \ldots, y_n)^\top$, $z=(z_1, \ldots, z_n)^\top\in S(D, D', n)$  for some~$n\in\N$ and $b\in\Z_{p}\setminus\set{0,-1}$ such that
$x + by + cz = 0$ with $c = -(b + 1)$. This implies that there exists some~$i$ with $1\leq i \leq n$ such that the component~$(x_i, y_i, z_i)$ of the vectors is a non-trivial weighted progression, i.e., it
is in~$P_{b}(D)$. However, the test above says that there is no triple
in~$P_{b}(D)$ which can occur, due to the fact that every digit in~$D'$ has to occur
$\abs{D}/n$ times in each vector. This is a contradiction to the assumption that the vectors are pairwise different. Thus, the set $S(D, D', n)$ is a cap for all suitable~$n\in\N$, and~$D$ is admissible.

For detailed examples we refer to Section~\ref{sec:proof}, cases $p = 11$ and $p = 17$.

\subsection{Matrix-Reducibility as a Sufficient Condition}

In order to show that a digit set~$D$ is admissible for some set~$D'\subseteq D$ of digits with fixed frequency, we can also use the following sufficient condition based on the matrix~$A$, which represents the linear constraints given in~\eqref{eq:equations-for-progressions} via $Ax = 0$. Again, we consider each equation $x + by + cz = 0$ separately and fix~$b\in\set{1,\ldots,p-2}$. Let $A_r$ be the reduced row echelon form of the matrix~$A$. For each row of~$A_r$ which only contains non-negative respectively non-positive entries, it is clear that the variables corresponding to non-zero entries of this row have to be zero (otherwise, the equation that corresponds to the said row cannot be fulfilled). This is due to the fact that we only search for non-negative solutions~$x$.

Thus, we can delete the columns of~$A_r$ that belong to these variables, and proceed with the next non-negative or non-positive row. Note that the deletion of columns can bring out new non-negative or non-positive rows. Naturally, this process determines if no such row is left in~$A_r$. If at the end all columns of~$A_r$ are deleted, then all variables~$x_i$ have to be zero. If this is the case for all $b\in\set{1,\ldots,p-2}$, then we say that~$(D, D')$ (or simply~$D$) is \emph{matrix-reducible}, which implies that the digit set~$D$ is admissible.

\begin{remark}
This procedure described here and the algorithm that we use for digit-reducibility in the previous subsection are \emph{essentially} of the same shape: While we start with the reduced row echelon form of~$A$ here, we can reformulate the algorithm of Subsection~\ref{sec:digit-red} in such a way that it is the same as this one but with~$A$ itself as initial matrix instead of its echelon form. The reason for the different descriptions of the two algorithms is our belief that it is easier and more convenient to handle with digits and weighted progressions instead of the corresponding matrices---at least if one wants to understand it and do it by hand.
\end{remark}

One can also think of other transformations of~$A$ as initial matrices for the reduction than the reduced row echelon form~$A_r$ or~$A$ itself, and even combine them. However, we refrained from optimizing this point because it works fine for our purpose.

For an example we refer to Section~\ref{sec:proof}, case $p = 23$.

\paragraph{} We remark that reducibility (both via digits or matrices) is only a sufficient condition for~$D$ to be
admissible, but not necessary.  The system of equations involved could have a
more sophisticated structure, and there are indeed admissible digit sets which are not
reducible. However, it turned out that these algorithmically simple tests are in fact very useful. They help to keep the proofs for the admissibility of digit sets simple and readable.

Moreover, digit- and matrix-reducibility are not equivalent: There exist digit sets which are digit- but not matrix-reducible (see case~$p = 17$ in Section~\ref{sec:proof}) and vice versa (see case~$p = 23$ in Section~\ref{sec:proof}).

Finally, it is of course also possible to combine the latter two approaches: We can choose between the digit- and matrix-reducibility algorithm depending on the parameter~$b$. Indeed, there are digit sets which are neither digit- nor matrix-reducible, but if we combine the two approaches, then the reducibility of the digit set can be shown.

\subsection{Elimination of Some Equations}

So far, it seems that admissibility (respectively reducibility) of a fixed digit set has to be checked in $p-2$ cases, namely for all equations~$x + by + cz \neq 0$ with $b$, $c\in\Z_p^n\setminus\set{0,-1}$ and $b + c = -1$. This is in fact not necessary: The following two observations help to \emph{significantly} reduce the cases that have to be studied later on.

\begin{remark}\label{rem:inverse-eq}
Let $p$ be a prime, $D\subseteq\Z_p$ and $b$, $c\in\Z_{p}$ with $c\neq 0$ and $b + c = -1$. Then the following assertions are true:
\begin{enumerate}[(a)]
\item\label{item:rem-a}
    A triple
  $(x,y,z)\in\Z_{p}^{3}$ satisfies $x + by + cz = 0$ if and only if
  $(z, y, x)$ satisfies $z + c^{-1}by + c^{-1}x = 0$. In particular,
  this means that $P_{b}(D)$ contains the same elements as
  $P_{c^{-1}b}(D)$ but mirror-inverted. Hence, only one of the two equations
  $x + by + cz = 0$ and $x + c^{-1}by + c^{-1}z = 0$ has to be considered.
\item\label{item:rem-b} The equation $x + by + cz = 0$ implies that $(x, y, z)\in P_b(D)$ holds if and only if $(x, z, y)\in P_c(D)$. In other words, $P_{b}(D)$ and
  $P_c(D)$ contain the same elements, but the last two components of the triples are always flipped. Thus, it is enough to consider one of the equations $x + by +cz = 0$ and $x + cy + bz = 0$. 
\end{enumerate}
\end{remark}

 We say that two equations $x + b_1 y + c_1 z = 0$ and $x + b_2 y + c_2 z = 0$ are equivalent if either $b_2 = c_1^{-1}b_1$ (case~\ref{item:rem-a} above) or $b_1 = c_2$ (case~\ref{item:rem-b} above).
Hence, only representatives of non-equivalent equations have to be tested for the cap set property. 

For the primes~$p$ considered in Theorem~\ref{thm:main}, the iterated application of the two cases of Remark~\ref{rem:inverse-eq} implies an immense simplification in our proof: It reduces the number of relevant equations from $p-2$ to $(p+1)/6$.

\section{Proof of Theorem~\ref{thm:main}}
\label{sec:proof}

If we find an admissible set of digits~$D$ and $D'\subseteq D$ of suitable sizes
(depending on~$p$), then the statements of the theorem follow by~\eqref{eq:cap-size}. Because of the comments above, it is enough to
show reducibility.

\begin{proof}[\normalfont\textbf{Case $p=11$}]
  We claim that $D = \set{0,1,3,4,5}$ with fixed digits $D' = \set{0,1,3}$ is digit-reducible (as well as matrix-reducible, which is not shown here), and study
  solutions $x$, $y$, $z\in D$ of $x+by+cz=0$ with
  $b\in \Z_{p}\setminus\set{0,-1}$ and $c = -(b + 1)$.

\begin{enumerate}
\item Case $x+z=2y$.  We list all triples of digits
  $(x,y,z)\in \{0,1,3,4,5\}^{3}$ that are solutions of $x+z=2y$, but leave out
  the trivial solutions $x=y=z$. These are the triples in $P_{-2}(D)$ and are
  given by
  \[(1,3,5),(3,4,5),(5,3,1),(5,4,3).\]

  We have $1\in D'$ and thus, the frequency of this digit has to be equal in any of the three positions. However,~$1$ does not occur in any of the triples in the second position, and as a consequence, the digits $1$ can only occur in the trivial progression~$(1,1,1)$. So the triples $(1,3,5)$ and $(5,3,1)$ cannot occur in any component of a potential weighted progression in~$S(D, D', n)$. 
  Hence, we delete $(1,3,5)$ and $(5,3,1)$ from the above list
  and
  \[(3,4,5),(5,4,3)\]
  remain.
  None of these two triples has the digit~$3$ in the second position. Thus, we delete both of them, and no triple from the set~$P_{-2}(D)$ remains.
  
  By Remark~\ref{rem:inverse-eq}~\ref{item:rem-b}, this also solves the case $x + 9z = 10y$. Moreover, as~$5$ is the inverse of~$9$ modulo~$11$, also the equation $x + 5z = 6y$ is covered due to Remark~\ref{rem:inverse-eq}~\ref{item:rem-a}.
  
\item Case $x+2z=3y$. For this equation ($b = -3$) the set of non-trivial
  weighted progressions~$P_{-3}(D)$ is given by
  \[(1, 0, 5),
 (1, 3, 4),
 (1, 4, 0),
 (3, 0, 4),
 (3, 1, 0),
 (4, 1, 5),
 (4, 5, 0),
 (5, 0, 3).\]

  As~$0$ never occurs in the first position and~$1$ never occurs in the last
  position, we can remove all triples with any occurrence of~$0$
  or~$1$. Therefore, again no non-trivial solutions in~$D$ remain.

  By Remark~\ref{rem:inverse-eq}~\ref{item:rem-a} with $2^{-1}\equiv 6 \bmod 11$, also the equation $x + 6z = 7y$ has no non-trivial solution in~$D$. By Remark~\ref{rem:inverse-eq}~\ref{item:rem-b}, this moreover solves the cases $x + 8z = 9y$ and $x + 4z = 5y$. Again applying the observation from Remark~\ref{rem:inverse-eq}~\ref{item:rem-a} to the latter two equations with $8^{-1}\equiv 7 \bmod 11$ respectively $4^{-1} \equiv 3 \bmod 11$, also the equations $x + 7z = 8y$ and $x + 3z = 4y$ are covered.
\end{enumerate}
Since we have (directly or via Remark~\ref{rem:inverse-eq}) considered all cases $b\in\set{1,\ldots,p-2}$, we conclude
that~$(D, D')$ is digit-reducible and thus, the appropriate size of~$S(D, D', n)$ follows by~\eqref{eq:cap-size}.
\end{proof}

\begin{proof}[\normalfont\textbf{Case $p=17$}]
  We claim that the digit set
  $D = \set{0, 1, 2, 4, 8, 9, 13}$ is reducible with fixed digits $D' = \set{0, 1, 2, 4, 8}$, and argue
  in analogy to the case~$p = 11$ above.
  
\begin{enumerate}
\item
Case $x+z=2y$.
We list all triples of digits in $P_{-2}(D)$, which are 
\begin{align*}
  &(0,1,2),(0,2,4),(0,4,8),(0,9,1),(0, 13, 9), (1, 9, 0), (1, 13, 8), (2, 1, 0),(4, 0, 13),\\
  &(4, 2, 0), (8, 0, 9), (8, 2, 13), (8, 4, 0), (8, 13, 1), (9, 0, 8), (9, 13, 0), (13, 0, 4), (13, 2, 8).
\end{align*}

Since the digit~$8$ does not occur in any of the triples on the second position, we can delete all triples that contain any~$8$ and obtain the remaining list
\begin{align*}
  &(0,1,2),(0,2,4),(0,9,1),(0, 13, 9), (1, 9, 0),\\ &(2, 1, 0),(4, 0, 13),
  (4, 2, 0),(9, 13, 0), (13, 0, 4).
\end{align*}

Next, we observe that no triple of this list has a~$4$ on the second
position. Thus, we delete all triples which contain the
digit~$4$. This yields
\begin{equation*}
  (0,1,2),(0,9,1),(0, 13, 9), (1, 9, 0),(2, 1, 0),(9, 13, 0).
\end{equation*}

Now this list contains no triple with the digits~$0$ or $2$ on the
second position. By deleting all triples containing these digits, no
non-trivial solution remains, which closes the argument for this case.

By Remark~\ref{rem:inverse-eq}, this also solves the cases 
\begin{itemize}
\item $x + 15z = 16y$ (as $2+ 16 \equiv 1\bmod 17$) and
\item $x+8z = 9y$ (as $15^{-1} \equiv 8\bmod 17$).
\end{itemize}

\item
Case $x+2z=3y$. This equation yields
\begin{align*}
  &(1,0,8), (1,9,13), (1,13,2), (2,1,9), (2,9,4), (4,1,8), (4,2,1),\\ &(4,13,9), (8,0,13), (8,4,2), (8,9,1), (9,0,4), (13,0,2),(13,4,8)
\end{align*}
as triples in $P_{-3}(D)$.
As $0$ never occurs in the first position and~$8$ never occurs in the
second position, we can remove all triples with any~$0$ or~$8$. The
remaining list is given by
\begin{equation*}
  (1,9,13), (1,13,2), (2,1,9), (2,9,4), (4,2,1),(4,13,9).
\end{equation*}
Next, we observe that the digit~$4$ never occurs in the second
position, which leads to the list 
\begin{equation*}
  (1,9,13), (1,13,2), (2,1,9).
\end{equation*}
Here, the digit $1$ does not occur in the third position. So all triples can be removed, which implies that there is no non-trivial
solution of~$x+2z=3y$ in~$D^{3}$.

Moreover, by repeatedly applying Remark~\ref{rem:inverse-eq}, this also solves the cases
\begin{multicols}{3}
\begin{itemize}
\item $x + 14z = 15y$,
\item $x + 9z = 10y$,
\item $x + 7z = 8y$,
\item $x + 11z = 12z$,
\item $x + 5z = 6y$.
\end{itemize}
\end{multicols}

\item Case $x + 3z = 4y$. This equation has the triples
  \begin{align*}
    &(1, 2, 8),
      (1, 13, 0),
      (2, 9, 0),
      (4, 1, 0),
      (8, 2, 0),
      (8, 13, 9),\\
    &(9, 1, 4),
      (9, 4, 8),
      (9, 8, 2),
      (13, 2, 4),
      (13, 4, 1),
      (13, 9, 2)
  \end{align*}
  as non-trivial solutions. Since~$0$ never occurs in the first
  position, we can remove
  all triples containing any~$0$ and obtain
  \begin{align*}
    (1, 2, 8),
      (8, 13, 9),
    (9, 1, 4),
      (9, 4, 8),
      (9, 8, 2),
      (13, 2, 4),
      (13, 4, 1),
      (13, 9, 2).
  \end{align*}
  Furthermore, the digits~$2$ and~$4$ do not occur in any triple in the first position,
  which leads to $(8, 13, 9)$ as the only remaining triple. 
  A single triple leads to the empty set.
  
  By Remark~\ref{rem:inverse-eq}, this also solves the cases
  \begin{multicols}{3}
\begin{itemize}
\item $x + 13z = 14y$,
\item $x + 6z = 7y$,
\item $x + 4z = 5y$,
\item $x + 10z = 11y$,
\item $x + 12z = 13y$.
\end{itemize}
\end{multicols}
      \end{enumerate}
      Since all cases $b\in\set{0,\ldots,p-2}$ are covered, this implies
that~$(D, D')$ is reducible and thus, also admissible. The appropriate size of the corresponding cap~$S(D, D', n)$ follows by~\eqref{eq:cap-size} again.
\end{proof}

\begin{proof}[\normalfont\textbf{Case $p=23$}]
    We claim that the digit set $D = \set{0, 1, 3, 4, 8, 9, 10, 12, 17}$ with fixed digits $D' = \set{0, 1, 3, 4, 8, 10, 17}$ is admissible. Unfortunately, $D$ is neither digit- nor matrix-reducible for any~$D'$ of size~$7$. So the admissibility has been checked by solving the corresponding IP with appropriate software, as described in Section~\ref{sec:admissible-sets}. This leads to the result
    \begin{equation*}
        C_{n,23} \gg \frac{9^n}{n^{3.5}},
    \end{equation*}
    as stated in Theorem~\ref{thm:main}.
    
    As a consolation prize, we show that $D$ is matrix-reducible for $D' = D$, i.e., if we fix the frequencies of \emph{all} digits. (This would lead to a lower bound of $9^n/n^4$.) For this purpose, we again study
    solutions $(x,y,z)\in D^3$ of $x+by+cz=0$ with
    $b\in \Z_{p}\setminus\set{0,-1}$ and $c = -(b + 1)$.
    
    The equivalent equations with respect to Remark~\ref{rem:inverse-eq} are given as follows, where each set represents an equivalence class:
    \begin{align*}
        &\set{x + z = 2y, x + 21z = 22y, x + 11z = 12y},\\
        &\set{x + 20z = 21y, x + 15z = 16y, x + 12z = 13y, x + 10z = 11y, x + 7z = 8y, x + 2z = 3y},\\
        &\set{x + 19z = 20y, x + 17z = 18y, x + 14z = 15y, x + 8z = 9y, x + 5z = 6y, x +3z = 4y},\\
        &\set{x + 18z = 19y, x + 16z = 6y, x + 13z = 14y, x + 9z = 10y, x + 6z = 7y, x + 4z = 5y}.
    \end{align*}
    Only one representative of each class has to be considered.
    
    Let us look at the equation $x + z = 2y$. Here, the progressions in~$P_{-2}(D)$ are given by
            \begin{align*}
                &(0, 4, 8),
 (0, 12, 1),
 (1, 9, 17),
 (1, 12, 0),
 (1, 17, 10),
 (3, 10, 17),
 (3, 17, 8),
 (4, 8, 12),\\
 &(8, 1, 17),
 (8, 4, 0),
 (8, 9, 10),
 (8, 10, 12),
 (8, 17, 3),
 (10, 9, 8),
 (10, 17, 1),\\
 &(12, 3, 17),
 (12, 8, 4),
 (12, 10, 8),
 (17, 1, 8),
 (17, 3, 12),
 (17, 9, 1),
 (17, 10, 3),
            \end{align*}
    and we call them $v_1$, \ldots, $v_{22}$ in the given order. The corresponding constraint matrix~$A$ (defined by the equations in~\eqref{eq:equations-for-progressions}) then has the form
        \begin{equation*} A = 
            \left(\begin{smallmatrix}
1 & 1 & 0 & 0 & 0 & 0 & 0 & 0 & 0 & 0 & 0 & 0 & 0 & 0 & 0 & 0 & 0 & 0 & 0 & 0 & 0 & 0 \\
0 & 0 & 1 & 1 & 1 & 0 & 0 & 0 & -1 & 0 & 0 & 0 & 0 & 0 & 0 & 0 & 0 & 0 & -1 & 0 & 0 & 0 \\
0 & 0 & 0 & 0 & 0 & 1 & 1 & 0 & 0 & 0 & 0 & 0 & 0 & 0 & 0 & -1 & 0 & 0 & 0 & -1 & 0 & 0 \\
-1 & 0 & 0 & 0 & 0 & 0 & 0 & 1 & 0 & -1 & 0 & 0 & 0 & 0 & 0 & 0 & 0 & 0 & 0 & 0 & 0 & 0 \\
0 & 0 & 0 & 0 & 0 & 0 & 0 & -1 & 1 & 1 & 1 & 1 & 1 & 0 & 0 & 0 & -1 & 0 & 0 & 0 & 0 & 0 \\
0 & 0 & -1 & 0 & 0 & 0 & 0 & 0 & 0 & 0 & -1 & 0 & 0 & -1 & 0 & 0 & 0 & 0 & 0 & 0 & -1 & 0 \\
0 & 0 & 0 & 0 & 0 & -1 & 0 & 0 & 0 & 0 & 0 & -1 & 0 & 1 & 1 & 0 & 0 & -1 & 0 & 0 & 0 & -1 \\
0 & -1 & 0 & -1 & 0 & 0 & 0 & 0 & 0 & 0 & 0 & 0 & 0 & 0 & 0 & 1 & 1 & 1 & 0 & 0 & 0 & 0 \\
0 & 0 & 0 & 0 & -1 & 0 & -1 & 0 & 0 & 0 & 0 & 0 & -1 & 0 & -1 & 0 & 0 & 0 & 1 & 1 & 1 & 1 \\
1 & 1 & 0 & -1 & 0 & 0 & 0 & 0 & 0 & -1 & 0 & 0 & 0 & 0 & 0 & 0 & 0 & 0 & 0 & 0 & 0 & 0 \\
0 & -1 & 1 & 1 & 1 & 0 & 0 & 0 & 0 & 0 & 0 & 0 & 0 & 0 & -1 & 0 & 0 & 0 & 0 & 0 & -1 & 0 \\
0 & 0 & 0 & 0 & 0 & 1 & 1 & 0 & 0 & 0 & 0 & 0 & -1 & 0 & 0 & 0 & 0 & 0 & 0 & 0 & 0 & -1 \\
0 & 0 & 0 & 0 & 0 & 0 & 0 & 1 & 0 & 0 & 0 & 0 & 0 & 0 & 0 & 0 & -1 & 0 & 0 & 0 & 0 & 0 \\
-1 & 0 & 0 & 0 & 0 & 0 & -1 & 0 & 1 & 1 & 1 & 1 & 1 & -1 & 0 & 0 & 0 & -1 & -1 & 0 & 0 & 0 \\
0 & 0 & 0 & 0 & 0 & 0 & 0 & 0 & 0 & 0 & 0 & 0 & 0 & 0 & 0 & 0 & 0 & 0 & 0 & 0 & 0 & 0 \\
0 & 0 & 0 & 0 & -1 & 0 & 0 & 0 & 0 & 0 & -1 & 0 & 0 & 1 & 1 & 0 & 0 & 0 & 0 & 0 & 0 & 0 \\
0 & 0 & 0 & 0 & 0 & 0 & 0 & -1 & 0 & 0 & 0 & -1 & 0 & 0 & 0 & 1 & 1 & 1 & 0 & -1 & 0 & 0 \\
0 & 0 & -1 & 0 & 0 & -1 & 0 & 0 & -1 & 0 & 0 & 0 & 0 & 0 & 0 & -1 & 0 & 0 & 1 & 1 & 1 & 1
\end{smallmatrix}\right),
        \end{equation*}
        where the first nine rows represent equations which arise from the first and second position in the vectors of~$P_{-2}(D)$ (left equation in~\eqref{eq:equations-for-progressions}), and the last nine rows represent the constraints for the positions one and three in the vectors of~$P_{-2}(D)$ (right equation in~\eqref{eq:equations-for-progressions}).
        
        Let us take a closer look at the construction of~$A$: For the first row of~$A$ we consider the first digit of~$D$, which is~$0$. This digit occurs in the triples $v_1 = (0,4,8)$ and $v_2 = (0, 12, 1)$ in the first position, and in none of the triples in the second position. Hence, following the left equation of~\eqref{eq:equations-for-progressions} with $d = 0$, this leads to the equation
        \begin{equation*}
            x_{v_1} + x_{v_2} = 0,
        \end{equation*}
        which is represented by the first row of~$A$.
        
        As a second, more sophisticated example, we consider the fourteenth row of~$A$ and the corresponding fifth digit in~$D$, which is~$8$. Now the first and the third positions of the progressions are significant (because the row is part of the last nine rows). In the vectors $v_9 = (8, 1, 17)$, $v_{10} = (8, 4, 0)$, $v_{11} = (8, 9, 10)$, $v_{12} = (8, 10, 12)$ and $v_{13} = (8, 17, 3)$, the digit~$8$ occurs in the first position. The vectors $v_1 = (0,4,8)$, $v_7 = (3, 17, 8)$, $v_{14} = (10,9,8)$, $v_{18} = (12,10,8)$ and $v_{19} = (17,1,8)$ contain~$8$ in the third position. Hence, due to the right equation in~\eqref{eq:equations-for-progressions} with $d = 8$, this yields the equation
        \begin{equation*}
            x_{v_{9}} + x_{v_{10}} + x_{v_{11}} + x_{v_{12}} + x_{v_{13}} = x_{v_{1}} + x_{v_{7}} + x_{v_{14}} + x_{v_{18}} + x_{v_{19}}.
        \end{equation*}
        This es exactly the equation represented by the fourteenth row of~$A$.
        
        The reduced row echelon form~$A_r$ of~$A$ is given by
                \begin{equation*} A_r = 
            \left(\begin{smallmatrix}
1 & 0 & 0 & 0 & 0 & 0 & 0 & 0 & 0 & 0 & 0 & 0 & 0 & 0 & -1 & 0 & 17 & 1 & 1 & 0 & 0 & 1 \\
0 & 1 & 0 & 0 & 0 & 0 & 0 & 0 & 0 & 0 & 0 & 0 & 0 & 0 & 1 & 0 & 6 & -1 & -1 & 0 & 0 & -1 \\
0 & 0 & 1 & 0 & 0 & 0 & 0 & 0 & 0 & 0 & 0 & 0 & 0 & 1 & 1 & 0 & 17 & -1 & -1 & 0 & 0 & -1 \\
0 & 0 & 0 & 1 & 0 & 0 & 0 & 0 & 0 & 0 & 0 & 0 & 0 & 0 & -1 & 0 & 18 & 1 & 1 & 0 & 0 & 1 \\
0 & 0 & 0 & 0 & 1 & 0 & 0 & 0 & 0 & 0 & 0 & 0 & 0 & -1 & 0 & 0 & 17 & -1 & -1 & 0 & -1 & -1 \\
0 & 0 & 0 & 0 & 0 & 1 & 0 & 0 & 0 & 0 & 0 & 0 & 0 & -1 & -1 & 0 & 21 & 1 & 0 & -1 & 0 & 1 \\
0 & 0 & 0 & 0 & 0 & 0 & 1 & 0 & 0 & 0 & 0 & 0 & 0 & 1 & 1 & 0 & 4 & 0 & 0 & 0 & 0 & -1 \\
0 & 0 & 0 & 0 & 0 & 0 & 0 & 1 & 0 & 0 & 0 & 0 & 0 & 0 & 0 & 0 & 22 & 0 & 0 & 0 & 0 & 0 \\
0 & 0 & 0 & 0 & 0 & 0 & 0 & 0 & 1 & 0 & 0 & 0 & 0 & 0 & 0 & 0 & 6 & -1 & 0 & 0 & -1 & -1 \\
0 & 0 & 0 & 0 & 0 & 0 & 0 & 0 & 0 & 1 & 0 & 0 & 0 & 0 & 1 & 0 & 5 & -1 & -1 & 0 & 0 & -1 \\
0 & 0 & 0 & 0 & 0 & 0 & 0 & 0 & 0 & 0 & 1 & 0 & 0 & 0 & -1 & 0 & 6 & 1 & 1 & 0 & 1 & 1 \\
0 & 0 & 0 & 0 & 0 & 0 & 0 & 0 & 0 & 0 & 0 & 1 & 0 & 0 & 0 & 0 & 2 & 0 & 0 & 1 & 0 & 0 \\
0 & 0 & 0 & 0 & 0 & 0 & 0 & 0 & 0 & 0 & 0 & 0 & 1 & 0 & 0 & 0 & 2 & 1 & 0 & -1 & 0 & 1 \\
0 & 0 & 0 & 0 & 0 & 0 & 0 & 0 & 0 & 0 & 0 & 0 & 0 & 0 & 0 & 1 & 2 & 1 & 0 & 0 & 0 & 0 \\
0 & 0 & 0 & 0 & 0 & 0 & 0 & 0 & 0 & 0 & 0 & 0 & 0 & 0 & 0 & 0 & 1 & 0 & 0 & 0 & 0 & 0 \\
0 & 0 & 0 & 0 & 0 & 0 & 0 & 0 & 0 & 0 & 0 & 0 & 0 & 0 & 0 & 0 & 0 & 0 & 0 & 0 & 0 & 0 \\
0 & 0 & 0 & 0 & 0 & 0 & 0 & 0 & 0 & 0 & 0 & 0 & 0 & 0 & 0 & 0 & 0 & 0 & 0 & 0 & 0 & 0 \\
0 & 0 & 0 & 0 & 0 & 0 & 0 & 0 & 0 & 0 & 0 & 0 & 0 & 0 & 0 & 0 & 0 & 0 & 0 & 0 & 0 & 0
\end{smallmatrix}\right).
        \end{equation*}
        Now we look for non-zero rows of~$A_r$ in which all entries are either non-negative or non-positive. The indices of these rows are given by $R_1 = \set{8,12,14,15}$ (and are framed in green in the following matrix). Next, we delete all columns of~$A_r$ in which any row~$r$ with $r\in R_1$ has some non-zero entry (symbolized by red lines). This means that we eliminate the variable which corresponds to this column. As a result, the first reduction step looks like
                        \begin{equation*}
            \left(\begin{smallmatrix}
1 & 0 & 0 & 0 & 0 & 0 & 0 & \tikzmark{o1}0 & 0 & 0 & 0 &\tikzmark{o2} 0 & 0 & 0 & -1 &\tikzmark{o4} 0 &\tikzmark{o5}17 &\tikzmark{o7} 1 & 1 &\tikzmark{o6} 0 & 0 & 1 \\
0 & 1 & 0 & 0 & 0 & 0 & 0 & 0 & 0 & 0 & 0 & 0 & 0 & 0 & 1 & 0 & 6 & -1 & -1 & 0 & 0 & -1 \\
0 & 0 & 1 & 0 & 0 & 0 & 0 & 0 & 0 & 0 & 0 & 0 & 0 & 1 & 1 & 0 & 17 & -1 & -1 & 0 & 0 & -1 \\
0 & 0 & 0 & 1 & 0 & 0 & 0 & 0 & 0 & 0 & 0 & 0 & 0 & 0 & -1 & 0 & 18 & 1 & 1 & 0 & 0 & 1 \\
0 & 0 & 0 & 0 & 1 & 0 & 0 & 0 & 0 & 0 & 0 & 0 & 0 & -1 & 0 & 0 & 17 & -1 & -1 & 0 & -1 & -1 \\
0 & 0 & 0 & 0 & 0 & 1 & 0 & 0 & 0 & 0 & 0 & 0 & 0 & -1 & -1 & 0 & 21 & 1 & 0 & -1 & 0 & 1 \\
0 & 0 & 0 & 0 & 0 & 0 & 1 & 0 & 0 & 0 & 0 & 0 & 0 & 1 & 1 & 0 & 4 & 0 & 0 & 0 & 0 & -1 \\
\tikzmark{l8}0 & 0 & 0 & 0 & 0 & 0 & 0 & 1 & 0 & 0 & 0 & 0 & 0 & 0 & 0 & 0 & 22 & 0 & 0 & 0 & 0 & 0 \tikzmark{r8}\\
0 & 0 & 0 & 0 & 0 & 0 & 0 & 0 & 1 & 0 & 0 & 0 & 0 & 0 & 0 & 0 & 6 & -1 & 0 & 0 & -1 & -1 \\
0 & 0 & 0 & 0 & 0 & 0 & 0 & 0 & 0 & 1 & 0 & 0 & 0 & 0 & 1 & 0 & 5 & -1 & -1 & 0 & 0 & -1 \\
0 & 0 & 0 & 0 & 0 & 0 & 0 & 0 & 0 & 0 & 1 & 0 & 0 & 0 & -1 & 0 & 6 & 1 & 1 & 0 & 1 & 1 \\
\tikzmark{l12}0 & 0 & 0 & 0 & 0 & 0 & 0 & 0 & 0 & 0 & 0 & 1 & 0 & 0 & 0 & 0 & 2 & 0 & 0 & 1 & 0 & 0 \tikzmark{r12}\\
0 & 0 & 0 & 0 & 0 & 0 & 0 & 0 & 0 & 0 & 0 & 0 & 1 & 0 & 0 & 0 & 2 & 1 & 0 & -1 & 0 & 1 \\
\tikzmark{l14}0 & 0 & 0 & 0 & 0 & 0 & 0 & 0 & 0 & 0 & 0 & 0 & 0 & 0 & 0 & 1 & 2 & 1 & 0 & 0 & 0 & 0\tikzmark{r14} \\
\tikzmark{l15}0 & 0 & 0 & 0 & 0 & 0 & 0 & 0 & 0 & 0 & 0 & 0 & 0 & 0 & 0 & 0 & 1 & 0 & 0 & 0 & 0 & 0 \tikzmark{r15}\\
0 & 0 & 0 & 0 & 0 & 0 & 0 & 0 & 0 & 0 & 0 & 0 & 0 & 0 & 0 & 0 & 0 & 0 & 0 & 0 & 0 & 0 \\
0 & 0 & 0 & 0 & 0 & 0 & 0 & 0 & 0 & 0 & 0 & 0 & 0 & 0 & 0 & 0 & 0 & 0 & 0 & 0 & 0 & 0 \\
0 & 0 & 0 & 0 & 0 & 0 & 0 & \tikzmark{u1}0 & 0 & 0 & 0 &\tikzmark{u2} 0 & 0 & 0 & 0 &\tikzmark{u4} 0 &\tikzmark{u5} 0 &\tikzmark{u7} 0 & 0 &\tikzmark{u6} 0 & 0 & 0
\end{smallmatrix}\right)
\DrawBox[ForestGreen, thick]{l8}{r8}
\DrawBox[ForestGreen, thick]{l12}{r12}
\DrawBox[ForestGreen, thick]{l14}{r14}
\DrawBox[ForestGreen, thick]{l15}{r15}
\DrawLine[Maroon, thick]{o1}{u1}
\DrawLine[Maroon, thick]{o2}{u2}
\DrawLine[Maroon, thick]{o4}{u4}
\DrawLineCorr[Maroon, thick]{o5}{u5}
\DrawLine[Maroon, thick]{o6}{u6}
\DrawLine[Maroon, thick]{o7}{u7}
\leadsto
    \left(\begin{smallmatrix}
1 & 0 & 0 & 0 & 0 & 0 & 0 & 0 & 0 & 0 & 0 & 0 & -1 & 1 & 0 & 1 \\
0 & 1 & 0 & 0 & 0 & 0 & 0 & 0 & 0 & 0 & 0 & 0 & 1 & -1 & 0 & -1 \\
0 & 0 & 1 & 0 & 0 & 0 & 0 & 0 & 0 & 0 & 0 & 1 & 1 & -1 & 0 & -1 \\
0 & 0 & 0 & 1 & 0 & 0 & 0 & 0 & 0 & 0 & 0 & 0 & -1 & 1 & 0 & 1 \\
0 & 0 & 0 & 0 & 1 & 0 & 0 & 0 & 0 & 0 & 0 & -1 & 0 & -1 & -1 & -1 \\
0 & 0 & 0 & 0 & 0 & 1 & 0 & 0 & 0 & 0 & 0 & -1 & -1 & 0 & 0 & 1 \\
0 & 0 & 0 & 0 & 0 & 0 & 1 & 0 & 0 & 0 & 0 & 1 & 1 & 0 & 0 & -1 \\
0 & 0 & 0 & 0 & 0 & 0 & 0 & 0 & 0 & 0 & 0 & 0 & 0 & 0 & 0 & 0 \\
0 & 0 & 0 & 0 & 0 & 0 & 0 & 1 & 0 & 0 & 0 & 0 & 0 & 0 & -1 & -1 \\
0 & 0 & 0 & 0 & 0 & 0 & 0 & 0 & 1 & 0 & 0 & 0 & 1 & -1 & 0 & -1 \\
0 & 0 & 0 & 0 & 0 & 0 & 0 & 0 & 0 & 1 & 0 & 0 & -1 & 1 & 1 & 1 \\
0 & 0 & 0 & 0 & 0 & 0 & 0 & 0 & 0 & 0 & 0 & 0 & 0 & 0 & 0 & 0 \\
0 & 0 & 0 & 0 & 0 & 0 & 0 & 0 & 0 & 0 & 1 & 0 & 0 & 0 & 0 & 1 \\
0 & 0 & 0 & 0 & 0 & 0 & 0 & 0 & 0 & 0 & 0 & 0 & 0 & 0 & 0 & 0 \\
0 & 0 & 0 & 0 & 0 & 0 & 0 & 0 & 0 & 0 & 0 & 0 & 0 & 0 & 0 & 0 \\
0 & 0 & 0 & 0 & 0 & 0 & 0 & 0 & 0 & 0 & 0 & 0 & 0 & 0 & 0 & 0 \\
0 & 0 & 0 & 0 & 0 & 0 & 0 & 0 & 0 & 0 & 0 & 0 & 0 & 0 & 0 & 0 \\
0 & 0 & 0 & 0 & 0 & 0 & 0 & 0 & 0 & 0 & 0 & 0 & 0 & 0 & 0 & 0
\end{smallmatrix}\right).
        \end{equation*}
        The rows~$r$ for $r\in R_1$ are zero now. Next, we proceed the same way with the resulting smaller matrix on the right-hand side. The only non-negative respectively non-positive non-zero row in this matrix is row~$13$. So we delete the corresponding two columns, and again proceed with the smaller matrix, and so on. The full remaining reduction is given by
        \allowdisplaybreaks
        \begin{multline*}
            \left(\begin{smallmatrix}
1 & 0 & 0 & 0 & 0 & 0 & 0 & 0 & 0 & 0 &\tikzmark{o21}0 & 0 & -1 & 1 & 0 &\tikzmark{o22}1 \\
0 & 1 & 0 & 0 & 0 & 0 & 0 & 0 & 0 & 0 & 0 & 0 & 1 & -1 & 0 & -1 \\
0 & 0 & 1 & 0 & 0 & 0 & 0 & 0 & 0 & 0 & 0 & 1 & 1 & -1 & 0 & -1 \\
0 & 0 & 0 & 1 & 0 & 0 & 0 & 0 & 0 & 0 & 0 & 0 & -1 & 1 & 0 & 1 \\
0 & 0 & 0 & 0 & 1 & 0 & 0 & 0 & 0 & 0 & 0 & -1 & 0 & -1 & -1 & -1 \\
0 & 0 & 0 & 0 & 0 & 1 & 0 & 0 & 0 & 0 & 0 & -1 & -1 & 0 & 0 & 1 \\
0 & 0 & 0 & 0 & 0 & 0 & 1 & 0 & 0 & 0 & 0 & 1 & 1 & 0 & 0 & -1 \\
0 & 0 & 0 & 0 & 0 & 0 & 0 & 0 & 0 & 0 & 0 & 0 & 0 & 0 & 0 & 0 \\
0 & 0 & 0 & 0 & 0 & 0 & 0 & 1 & 0 & 0 & 0 & 0 & 0 & 0 & -1 & -1 \\
0 & 0 & 0 & 0 & 0 & 0 & 0 & 0 & 1 & 0 & 0 & 0 & 1 & -1 & 0 & -1 \\
0 & 0 & 0 & 0 & 0 & 0 & 0 & 0 & 0 & 1 & 0 & 0 & -1 & 1 & 1 & 1 \\
0 & 0 & 0 & 0 & 0 & 0 & 0 & 0 & 0 & 0 & 0 & 0 & 0 & 0 & 0 & 0 \\
\tikzmark{l21}0 & 0 & 0 & 0 & 0 & 0 & 0 & 0 & 0 & 0 & 1 & 0 & 0 & 0 & 0 & 1\tikzmark{r21} \\
0 & 0 & 0 & 0 & 0 & 0 & 0 & 0 & 0 & 0 & 0 & 0 & 0 & 0 & 0 & 0 \\
0 & 0 & 0 & 0 & 0 & 0 & 0 & 0 & 0 & 0 & 0 & 0 & 0 & 0 & 0 & 0 \\
0 & 0 & 0 & 0 & 0 & 0 & 0 & 0 & 0 & 0 & 0 & 0 & 0 & 0 & 0 & 0 \\
0 & 0 & 0 & 0 & 0 & 0 & 0 & 0 & 0 & 0 & 0 & 0 & 0 & 0 & 0 & 0 \\
0 & 0 & 0 & 0 & 0 & 0 & 0 & 0 & 0 & 0 &\tikzmark{u21}0 & 0 & 0 & 0 & 0 &\tikzmark{u22}0
\end{smallmatrix}\right)
\DrawBox[ForestGreen, thick]{l21}{r21}
\DrawLine[Maroon, thick]{o21}{u21}
\DrawLine[Maroon, thick]{o22}{u22}
\leadsto
            \left(\begin{smallmatrix}
1 & 0 & 0 & 0 & 0 & 0 &\tikzmark{o31} 0 & 0 & 0 & 0 &\tikzmark{o33} 0 &\tikzmark{o32} -1 & 1 & 0 \\
0 & 1 & 0 & 0 & 0 & 0 & 0 & 0 & 0 & 0 & 0 & 1 & -1 & 0 \\
0 & 0 & 1 & 0 & 0 & 0 & 0 & 0 & 0 & 0 & 1 & 1 & -1 & 0 \\
0 & 0 & 0 & 1 & 0 & 0 & 0 & 0 & 0 & 0 & 0 & -1 & 1 & 0 \\
0 & 0 & 0 & 0 & 1 & 0 & 0 & 0 & 0 & 0 & -1 & 0 & -1 & -1 \\
0 & 0 & 0 & 0 & 0 & 1 & 0 & 0 & 0 & 0 & -1 & -1 & 0 & 0 \\
\tikzmark{l31}0 & 0 & 0 & 0 & 0 & 0 & 1 & 0 & 0 & 0 & 1 & 1 & 0 & 0\tikzmark{r31} \\
0 & 0 & 0 & 0 & 0 & 0 & 0 & 0 & 0 & 0 & 0 & 0 & 0 & 0 \\
0 & 0 & 0 & 0 & 0 & 0 & 0 & 1 & 0 & 0 & 0 & 0 & 0 & -1 \\
0 & 0 & 0 & 0 & 0 & 0 & 0 & 0 & 1 & 0 & 0 & 1 & -1 & 0 \\
0 & 0 & 0 & 0 & 0 & 0 & 0 & 0 & 0 & 1 & 0 & -1 & 1 & 1 \\
0 & 0 & 0 & 0 & 0 & 0 & 0 & 0 & 0 & 0 & 0 & 0 & 0 & 0 \\
0 & 0 & 0 & 0 & 0 & 0 & 0 & 0 & 0 & 0 & 0 & 0 & 0 & 0 \\
0 & 0 & 0 & 0 & 0 & 0 & 0 & 0 & 0 & 0 & 0 & 0 & 0 & 0 \\
0 & 0 & 0 & 0 & 0 & 0 & 0 & 0 & 0 & 0 & 0 & 0 & 0 & 0 \\
0 & 0 & 0 & 0 & 0 & 0 & 0 & 0 & 0 & 0 & 0 & 0 & 0 & 0 \\
0 & 0 & 0 & 0 & 0 & 0 & 0 & 0 & 0 & 0 & 0 & 0 & 0 & 0 \\
0 & 0 & 0 & 0 & 0 & 0 &\tikzmark{u31} 0 & 0 & 0 & 0 &\tikzmark{u33} 0 &\tikzmark{u32} 0 & 0 & 0
\end{smallmatrix}\right)
    \DrawBox[ForestGreen, thick]{l31}{r31}
    \DrawLine[Maroon, thick]{o31}{u31}
    \DrawLineCorrTwo[Maroon, thick]{o32}{u32}
    \DrawLine[Maroon, thick]{o33}{u33}
    \leadsto\\
                \leadsto
                \left(\begin{smallmatrix}
\tikzmark{l41}1 & 0 & 0 &\tikzmark{o42} 0 & 0 &\tikzmark{o43} 0 & 0 & 0 &\tikzmark{o44} 0 &\tikzmark{o45} 1 &\tikzmark{o46} 0\tikzmark{r41} \\
0 & 1 & 0 & 0 & 0 & 0 & 0 & 0 & 0 & -1 & 0 \\
0 & 0 & 1 & 0 & 0 & 0 & 0 & 0 & 0 & -1 & 0 \\
\tikzmark{l42}0 & 0 & 0 & 1 & 0 & 0 & 0 & 0 & 0 & 1 & 0\tikzmark{r42} \\
0 & 0 & 0 & 0 & 1 & 0 & 0 & 0 & 0 & -1 & -1 \\
\tikzmark{l43}0 & 0 & 0 & 0 & 0 & 1 & 0 & 0 & 0 & 0 & 0\tikzmark{r43} \\
0 & 0 & 0 & 0 & 0 & 0 & 0 & 0 & 0 & 0 & 0 \\
0 & 0 & 0 & 0 & 0 & 0 & 0 & 0 & 0 & 0 & 0 \\
0 & 0 & 0 & 0 & 0 & 0 & 1 & 0 & 0 & 0 & -1 \\
0 & 0 & 0 & 0 & 0 & 0 & 0 & 1 & 0 & -1 & 0 \\
\tikzmark{l44}0 & 0 & 0 & 0 & 0 & 0 & 0 & 0 & 1 & 1 & 1\tikzmark{r44} \\
0 & 0 & 0 & 0 & 0 & 0 & 0 & 0 & 0 & 0 & 0 \\
0 & 0 & 0 & 0 & 0 & 0 & 0 & 0 & 0 & 0 & 0 \\
0 & 0 & 0 & 0 & 0 & 0 & 0 & 0 & 0 & 0 & 0 \\
0 & 0 & 0 & 0 & 0 & 0 & 0 & 0 & 0 & 0 & 0 \\
0 & 0 & 0 & 0 & 0 & 0 & 0 & 0 & 0 & 0 & 0 \\
0 & 0 & 0 & 0 & 0 & 0 & 0 & 0 & 0 & 0 & 0 \\
\tikzmark{u41}0 & 0 & 0 &\tikzmark{u42} 0 & 0 &\tikzmark{u43} 0 & 0 & 0 &\tikzmark{u44} 0 &\tikzmark{u45} 0 &\tikzmark{u46} 0
\end{smallmatrix}\right)
\DrawBox[ForestGreen, thick]{l41}{r41}
\DrawBox[ForestGreen, thick]{l42}{r42}
\DrawBox[ForestGreen, thick]{l43}{r43}
\DrawBox[ForestGreen, thick]{l44}{r44}
\DrawLine[Maroon, thick]{l41}{u41}
\DrawLine[Maroon, thick]{o42}{u42}
\DrawLine[Maroon, thick]{o43}{u43}
\DrawLine[Maroon, thick]{o44}{u44}
\DrawLine[Maroon, thick]{o45}{u45}
\DrawLine[Maroon, thick]{o46}{u46}
\leadsto
\left(\begin{smallmatrix}
0 & 0 & 0 & 0 & 0 \\
1 & 0 & 0 & 0 & 0 \\
0 & 1 & 0 & 0 & 0 \\
0 & 0 & 0 & 0 & 0 \\
0 & 0 & 1 & 0 & 0 \\
0 & 0 & 0 & 0 & 0 \\
0 & 0 & 0 & 0 & 0 \\
0 & 0 & 0 & 0 & 0 \\
0 & 0 & 0 & 1 & 0 \\
0 & 0 & 0 & 0 & 1 \\
0 & 0 & 0 & 0 & 0 \\
0 & 0 & 0 & 0 & 0 \\
0 & 0 & 0 & 0 & 0 \\
0 & 0 & 0 & 0 & 0 \\
0 & 0 & 0 & 0 & 0 \\
0 & 0 & 0 & 0 & 0 \\
0 & 0 & 0 & 0 & 0 \\
0 & 0 & 0 & 0 & 0
\end{smallmatrix}\right)
\leadsto
\left(\begin{smallmatrix}
\ \\
\ \\
\ \\
\ \\
\ \\
\ \\
\ \\
\ \\
\ \\
\ \\
\ \\
\ \\
\ \\
\ \\
\ \\
\ \\
\ \\
\ \\
\ \\
\
\end{smallmatrix}\right),
        \end{multline*}
        where the last step trivially follows.
        
        All other equations can be handled analogously.
\end{proof}

\begin{proof}[\normalfont\textbf{Case $p \in \set{29, 41}$}]
  For these moduli, we only give the digit sets for which the
  digit-reducibility can be checked analogously to the cases $p=11$ and $p=17$: The sets $(D, D')$ with
  \begin{itemize}
      \item $D = \set{0, 1, 2, 3, 4, 6, 14, 16, 22, 26}$ with $D' = \set{1, 2, 3, 4, 6, 16, 22, 26}$ for $p = 29$ and
      \item $D = \set{1, 2, 4, 5, 6, 9, 15, 16, 27, 32, 33, 35}$ with $D' = \set{1, 2, 4, 5, 6, 9, 15, 27, 32, 33}$ for $p = 41$
  \end{itemize}
  are digit-reducible.
  
  Finally, we give the lists of equivalent equations with respect to Remark~\ref{rem:inverse-eq}; each set of equations represents an equivalence class. In the case $p = 29$, we have
  \begin{align*}
      &\set{x + z = 2y, x + 27z = 28y, x + 14z = 15y},\\
      &\set{x+26z=27y,
x+19z=20y,
x+15z=16y,
x+13z=14y,
x+9z=10y,
x+2z=3y
},\\
      &\set{x+25z=26y,
x+21z=22y,
x+18z=19y,
x+10z=11y,
x+7z=8y,
x+3z=4y},\\
&\set{x+24z=25y,
x+23z=24y,
x+22z=23y,
x+6z=7y,
x+5z=6y,
x+4z=5y},\\
&\set{x+20z=21y,
x+17z=18y,
x+16z=17y,
x+12z=13y,
x+11z=12y,
x+8z=9y},
  \end{align*}
  and for $p = 41$, the classes are given by
  \begin{align*}
      &\set{x+z=2y, x+39z=40y,
x+20z=21y},\\
      &\set{x+38z=39y,
x+27z=28y,
x+21z=22y,
x+19z=20y,
x+13z=14y,
x+2z=3y},\\
      &\set{x+37z=38y,
x+30z=31y,
x+26z=27y,
x+14z=15y,
x+10z=11y,
x+3z=4y},\\
      &\set{x+36z=37y,
x+32z=33y,
x+31z=32y,
x+9z=10y,
x+8z=9y,
x+4z=5y},\\
      &\set{x+35z=36y,
x+34z=35y,
x+33z=34y,
x+7z=8y,
x+6z=7y,
x+5z=6y},\\
      &\set{x+29z=30y,
x+25z=26y,
x+23z=24y,
x+17z=18y,
x+15z=16y,
x+11z=12y},\\
      &\set{x+28z=29y,
x+24z=25y,
x+22z=23y,
x+18z=19y,
x+16z=17y,
x+12z=13y}.
  \end{align*}
  Only one representative of each of these classes has to be considered.

  This concludes the proof of Theorem~\ref{thm:main}.
\end{proof}

\section{Non-Equivalent Caps}{\label{sec:non-equivalent}}
Two caps are equivalent if there is an affine transformation from one cap to the other. In some cases, two caps in $\ag(n,p)$ based on the above digit constructions but with different digit sets $D_1$ and $D_2$ are equivalent, while in other cases they are not. We briefly discuss this in the cases $p=5$ and $p=11$.

If a digit set $D_1\subseteq \Z_p$ can be mapped by an affine transformation $f(x)=ax+b$ to another digit set $D_2\subseteq \Z_p$, then the corresponding caps are equivalent.

For example, modulo $p=5$ all digit sets consisting of three distinct digits are equivalent. One can first map two arbitrary digits to $0$ and $1$.
Then the three remaining digit sets
$D_1=\{0,1,2\}$, $D_2=\{0,1,3\}$ and $D_3=\{0,1,4\}$ can be seen to be equivalent:
$D_3$ is mapped to $D_1$ by $f(x)=x+1$ and
$D_2$ is mapped to $D_1$ by $f(x)=3x+2$.

A simple criterion to see that two digit sets are not equivalent is as follows:
For a given digit set write the multiset of differences (including the gap from the largest digit to $p$).
If the multiset of differences of two digit sets $D_1$ and $D_2$ contain different frequencies of differences, then the two digit sets are not equivalent.

Applying this modulo~$5$ to the above digit sets gives twice the set of differences $\{1,1,3\}$ and once $\{1,2,2\}$. This helps finding the map $f(x)=3x+2$.

On the other hand, we easily find many admissible digit sets are not equivalent modulo~$11$:
$D_1=\{0, 1, 2, 3, 4\}$  with difference multiset $\{1,1,1,1,7\}$
is non-equivalent to
$D_2=\{0, 1, 2, 3, 6\}$ with difference multiset $\{1,1,1,3,5\}$,
and both are different from $D_3=\{0, 1, 2, 3, 7\}$ with difference multiset $\{1,1,1,4,4\}$.

Another criterion is that the order of gaps of the same frequencies must also be preserved:
$D_4= \{0, 1, 2, 6, 7\}$ is different from the earlier three digit sets,
as $D_4$ does not contain four elements in an arithmetic progression, which would be preserved by an affine map.
We leave it to the reader to argue why $D_5=\{0, 1, 2, 6, 8\}$
and $D_6=\{0, 1, 2, 8, 9\}$ lead to further non-equivalent digit sets.

For larger primes, the number of admissible digit sets is typically much larger than the number $p(p-1)$ of affine transformations of the digit sets.
Hence, our digit-based constructions typically indicate the existence of many non-equivalent caps with the same number of points. (However, we do not formally prove these caps are non-equivalent.) In any case this seems to be of interest even for those primes for which these caps are not larger than previously known ones.

\textbf{Postscriptum:} In the meantime, we have improved some of the results of this paper. For example, we have constructed caps in $\F_5^n$ with asymptotic growth at least $3.23^n$. In this case, the corresponding $\mu$ is actually larger than the current record for $p = 3$; see~\cite{edel:2004:product-caps, tyrell}. 

\bibliography{caps}
\bibliographystyle{amsplainurl}

\end{document}